\documentclass[12pt, twoside, leqno]{article}
\usepackage{amsmath,amsthm}
\usepackage{amssymb,latexsym}
\usepackage{enumerate}

\newtheorem{thm}{Theorem}[section]
\newtheorem{prop}[thm]{Proposition}
\newtheorem{cor}[thm]{Corollary}
\newtheorem{lem}[thm]{Lemma}
\newtheorem{prob}[thm]{Problem}
\newtheorem{ques}[thm]{Question}

\theoremstyle{definition}
\newtheorem{defin}[thm]{Definition}
\newtheorem{rem}[thm]{Remark}
\newtheorem{exa}[thm]{Example}

\numberwithin{equation}{section}

\usepackage{mathrsfs}
\usepackage{CJK}
\usepackage{amsmath}
\usepackage{fancyhdr}

\newcommand{\rank}{\operatorname{rank}}
\newcommand{\corank}{\operatorname{corank}}
\newcommand{\diag}{\operatorname{diag}}

\begin{document}

\baselineskip=17pt

\title{Some Remarks on the Jacobian Conjecture and Dru{\.z}kowski mappings}
\author{Dan Yan \\
School of Mathematical Sciences, Graduate University of \\
Chinese Academy of Sciences, Beijing 100049, China \\
\emph{E-mail:} yan-dan-hi@163.com \\
\and
Michiel de Bondt\footnote{The second author was supported by the Netherlands
                          Organisation for Scientific Research (NWO).} \\
Department of Mathematics, Radboud University \\
Nijmegen, The Netherlands \\
\emph{E-mail:} M.deBondt@math.ru.nl}
\date{}

\maketitle

\renewcommand{\thefootnote}{}

\renewcommand{\thefootnote}{\arabic{footnote}}
\setcounter{footnote}{0}

\begin{abstract} In this paper, we first show that the Jacobian Conjecture is true for non-homogeneous
power linear mappings under some conditions. Secondly, we prove an equivalent statement about the Jacobian
Conjecture in dimension $r\geq 1$ and give some partial results for $r=2$.

Finally, for a homogeneous power linear Keller map $F=X+H$ of degree $d \ge 2$,
we give the inverse polynomial map under the condition that $JH^3=0$.
We shall show that $\deg (F^{-1})\leq d^k$ if $k \le 2$ and $JH^{k+1}=0$,
but also give an example with $d = 2$ and $JH^4=0$ such that $\deg (F^{-1})> d^3$.
\end{abstract}
 \vskip 2.5mm

\section{Introduction}

Throughout this paper, we will write {\bf K} for any field of characteristic zero and
${\bf K}[x]={\bf K}[x_1,x_2,\ldots,x_r]$ (${\bf K}[X]={\bf K}[x_1,x_2,\ldots,x_n]$)
for the polynomial algebra over {\bf K} with $r$ ($n$) indeterminates $x=x_1,x_2,\ldots,x_r(,\ldots x_n)$.
Let $f=(f_1,f_2,\ldots,\allowbreak f_r):{\bf{K}}^r\rightarrow{\bf{K}}^r$
($F=(F_1,F_2,\ldots,F_n):{\bf{K}}^n\rightarrow{\bf{K}}^n)$
be a polynomial mapping, that is, $f_i\in{\bf{K}}[x]$ for all $1\leq i\leq r$
($F_i\in{\bf{K}}[X]$ for all $1\leq i\leq n$).
Let $Jf=(\frac{\partial f_i}{\partial x_j})_{r\times r}$ and
let $JF=(\frac{\partial F_i}{\partial x_j})_{n\times n}$ be the
Jacobian matrix of $F$. Write $M|_{G}$ for replacing either $x$ or $X$ by $G$
in $M$ (only one of them will possible if $x \ne X$).

The Jacobian Conjecture (JC)
raised by O.H. Keller in 1939 in \cite{MR1550818} states that a polynomial
mapping $F:{\bf{K}}^n\rightarrow{\bf{K}}^n$ is invertible if the
Jacobian determinant $\det JF$ is a nonzero constant. This conjecture
has been attacked by many people from various research fields, but it is still open, even for $n\geq 2$.
Only the case $n=1$ is obvious. For more information about the wonderful 70-year history, see \cite{MR0663785},
\cite{MR1790619}, and the references therein. It can easily be seen that the JC
is true if the JC holds for all polynomial mappings whose Jacobian
determinant is 1. We make use of this convention in the present paper.

In 1980, S.S.S.Wang (\cite{MR0585736}) showed that
the JC holds for all polynomial mappings of degree 2 in all dimensions.
The most powerful result is the reduction to degree
3, due to H.Bass, E.Connell and D.Wright (\cite{MR0663785}) in 1982 and
A.Yagzhev (\cite{MR0592226}) in 1980, which asserts that the JC is true if the JC
holds for all polynomial mappings of degree 3 (what is more, if the JC
holds for all cubic homogeneous polynomial mappings!). It is even shown in \cite{MR2138860} that the condition that
$JH$ is symmetric and $H$ is cubic homogeneous is sufficient.
In the same spirit of the above degree reduction method, another efficient way to tackle the JC is the
Dru{\.z}kowski's Theorem (\cite{MR0714105}): the JC is true if it is true for all Dru{\.z}kowski mappings
(in all dimension $\geq 2$).
One more interesting result is due to Gorni-Zampieri (\cite{MR1621913}), who proved in 1997 that there exist
Gorni-Zampieri pairings between the cubic homogeneous polynomial mappings and the Dru{\.z}kowski mappings.

Recall that $F$ is a cubic homogeneous mapping if $F=X+H$ with
$X$ the identity (written as a column vector) and each component of
$H$ being either zero or cubic homogeneous. A cubic homogeneous
mapping $F=X+H$ is a {\it\bf Dru{\.z}kowski (or cubic linear) mapping}
if each component of $H$ is either zero or a third power of a linear
form. Each Dru{\.z}kowski mapping $F$ is associated to a scalar matrix
$A$ such that $F=X+(AX)^{*3}$, where $(AX)^{*3}$ is the {\it\bf
Dru{\.z}kowski symbol} for the vector $((A_1X)^3,\ldots, (A_nX)^3)$ with
$A_i$ the $i$-th row of $A$. Clearly, a Dru{\.z}kowski mapping is
uniquely determined by this matrix $A$. In section 2, we prove that the JC is true for Dru{\.z}kowski
mappings in some cases.

Apparently, the notion of a Dru{\.z}kowski mapping can be easily generalized. Namely, for any
positive integer $d\geq 2$, we say that $F=X+H$ is homogeneous power linear of degree $d$
if each component of $H$ is either zero or a $d$-th power of a linear form. The JC is true in general
if it is true for homogeneous power linear maps of degree $d$, where $d$ is any integer larger than two.
If $F$ is an invertible polynomial map of degree $d$ in dimension $n$, the degree of its inverse is
at most $d^{n-1}$.
This has been proved in \cite{MR0663785}. But if additionally $F=X+H$ is homogeneous power linear of degree $d$
such that $JH^3=0$, then the degree of the inverse of $F$ is at most $d^2$. We will prove this in section 5,
using results of section 3.

In section 3, we generalize the definition of GZ-paired in \cite{homokema}, \cite{MR1790619} and \cite{MR1621913}.
We use this to prove in section 4 that
the Jacobian Conjecture in dimension $r\geq 1$ is equivalent to the Jacobian Conjecture for
non-homogeneous power linear maps with $\rank A\leq r$ and prove the Jacobian Conjecture
is true in this case for $r=2$ under the condition that $\det (DJH+I)=1$, where
$D$ is a certain diagonal matrix.

\section{The JC for Dru{\.z}kowski mappings}

\begin{thm}
Let $F=X+H$ such that $H_i=(A_iX)^{d_i}$ is a power of a linear form for each $i$, where $A_i$ is the
$i$-th row of a matrix $A$. If ${\operatorname{Tr}}JH=0$ and all the determinants of the
$i\times i$ principal minors of $A$ are zero for $2\leq i\leq n$, then $F$ is a polynomial automorphism.
\end{thm}

\begin{proof}
Since $JH(y^{(1)})+\cdots+JH(y^{(n)})=DA$, where $D$ is some diagonal matrix and
$y^{(i)}\in {\bf K}^n$ for $1\leq i\leq n$, we have that all the determinants of the
$i\times i$ principal minors of $JH(y^{(1)})+\cdots+JH(y^{(n)})$ are zero for $2\leq i\leq n$.
Furthermore, the trace of $JH(y^{(1)})+\cdots+JH(y^{(n)})$ is zero by additivity.
Hence $JH(y^{(1)})+\cdots+JH(y^{(n)})$ is nilpotent and
$\det (JF(y^{(1)})+\cdots+JF(y^{(n)}))=
\det (nI_n+JH(y^{(1)})+\cdots+JH(y^{(n)}))=n^n$.
Thus we deduce from \cite{MR2948624} (Theorem 3.5) that $F$ is invertible.
\end{proof}

\begin{cor}
Let $F=X+H$ be a Dru{\.z}kowski mapping, say that $H_i = (A_iX)^3$ for each $i$.
If $\det JF=1$ and all the determinants of the $i\times i$ principal minors of $A$
are zero for $2\leq i\leq n-4$, then F is a polynomial automorphism.
\end{cor}
\begin{proof}
If there exists an $i \in \{n-3,n-2,n-1,n\}$ such that some
$i\times i$ principal minor of $A$ is nonzero, then
$\corank A\leq 3$. Therefore, $F$ is a tame automorphism in that case, see
\cite{homokema} (Theorem 7.1.1).

Since $\det JF=1$, we have ${\operatorname{Tr}}JH=0$. If all the determinants of the
$i\times i$ principal minors of $A$ are zero for $ i\geq n-3$, then the conclusion follows from Theorem 2.1.
\end{proof}

\begin{cor}
Let $F=X+H$ be a Dru{\.z}kowski mapping in dimension $n$. If $\det JF=1$ and the diagonal
of $JH$ is entirely nonzero, then $F$ is tame for $n\leq 9$ and linearly triangularizable for $n\leq 7$.
\end{cor}
\begin{proof}
Since $\det JF=1$ and the diagonal of $JH$ is entirely
nonzero, we have $\rank A\leq [\frac{n}{2}]$, see \cite{MR2810555} (Theorem).
Thus $\rank A\leq 4$ when $n\leq 9$, in which case $F$ is tame, see \cite{homokema}
(Theorem 7.1.2). Furthermore, $\rank A\leq 3$ when $n\leq 7$,
in which case $F$ is linearly triangularizable, see \cite{MR2179727} (Corollary 4.1).
\end{proof}

\section{Gorni-Zampieri pairing}

In the rest of this paper, $\ker_{\bf L} M$ will be a vector space over the field ${\bf L}$ with
coordinates in ${\bf L}$ (where $M$ is a matrix over ${\bf L}$), and $\dim_{\bf L} V$ will
be the dimension of $V$ over ${\bf L}$ (where $V$ is a vector space over ${\bf L}$).
If the subscript field ${\bf L}$ is omitted, then ${\bf L} = {\bf K}(x)$.
So $\ker_{\bf K} M = (\ker M) \cap {\bf K}^n$ if $M$ has $n$ columns.

\begin{defin} \label{3.1}
Let $f: {\bf K}^r\rightarrow {\bf K}^r$ be polynomial maps and
$F: {\bf K}^n\rightarrow {\bf K}^n$ be nonhomogeneous power-linear maps with $n>r$.
We say that $f$ and $F$ are GZ-paired (weakly GZ-paired) through
the matrices $B\in M_{r,n}({\bf K})$ and $C\in M_{n,r}({\bf K})$ if
\begin{enumerate}
\item[1)] $f(y)=BF(Cy)$ for all $y\in {\bf K}^r$,
\item[2)] $BC=I_r$,
\item[3)] $\ker B=\ker JH$
($\ker B\subseteq \ker JH$),
\end{enumerate}
where $H = F - X$.
\end{defin}

From lemma \ref{3.2} below, we deduce that definition \ref{3.1} above is a generalization of
\cite{MR1790619} (Definition 6.4.1) and \cite{MR1621913} (Definition 1.2).

\begin{lem} \label{3.2}
Let $H=(H_1,H_2,\ldots,H_n)^t$ and $H_i$ is a power of $A_iX$ for each $i$,
where $A_i$ is the $i$-th row of $A$. Then $\ker JH=\ker A$ and
$n-\dim_{\bf K} (\ker_{\bf K} JH)=\rank A$.
\end{lem}
\begin{proof}
Since $JH=\diag (d_1t_1^{d_1-1},\ldots,d_nt_n^{d_n-1})A$, where $t_i=A_iX$ for $1\leq i\leq n$,
we have $\ker JH=\ker A$ and $n-\dim_{\bf K} (\ker_{\bf K} JH)
= n-\dim_{\bf K} ((\ker JH) \cap {\bf K}^n) \ge n - \dim (\ker JH) = \rank A$.
\end{proof}

\begin{thm} \label{3.3}
Let $f: {\bf K}^r\rightarrow {\bf K}^r$ be a (non)homogeneous polynomial map of degree (at most) $d$.
Then there exists an $n > r$ and a (non)homogeneous power linear map $F$ of degree $d$ such that
$f$ and $F$ are GZ-paired through some matrices $B\in M_{r,n}({\bf K})$ and $C\in M_{n,r}({\bf K})$.
\end{thm}
\begin{proof}
The proof is similar to the homogeneous case of \cite{homokema} (Theorem 6.2.8),
and the cubic homogeneous case of \cite{MR1790619} (Theorem 6.4.2) and \cite{MR1621913}
(Theorem 1.3).
\end{proof}

\begin{thm} \label{3.4}
Let $F: {\bf K}^n\rightarrow {\bf K}^n$ be a (non)homogeneous polynomial map of degree $d$ and let
$r\geq n-\dim_{\bf K} (\ker_{\bf K} JH)$, where $H=F-X$. If $0<r<n$, then there
exists a (not necessarily) homogeneous polynomial map $f$ of degree (at most) $d$ in dimension $r$,
such that $f$ and $F$ are weakly GZ-paired through
some matrices $B\in M_{r,n}({\bf K})$ and $C\in M_{n,r}({\bf K})$.
\end{thm}
\begin{proof}
Since $\dim (\ker_{\bf K} JH)\geq n-r$,
$\ker_{\bf K} JH$ has a linear subspace $S$ of dimension $n-r$.
Take $B$ in $M_{r,n}({\bf K})$ such that $\ker_{\bf K} B=S$.
Then $\rank B=r$. From $\ker_{\bf K} B \subseteq \ker_{\bf K} JH$ and
the fact that $\ker B$ is generated by vectors over $K$,
we obtain $\ker B \subseteq \ker JH$. Since $\rank B=r$,
there exists a $C$ in $M_{r,n}({\bf K})$ such that $BC=I_r$. Now $f:=BF(Cx)$ has the desired properties.
\end{proof}

\begin{cor} \label{3.5}
Let $F:{\bf K}^n\rightarrow {\bf K}^n$ be a (non)homogeneous power linear map of degree $d$ and let
$r=\rank A$, where $A$ is defined by $F_i-X_i=(A_iX)^{d_i}$. If $F$ is of Keller type
and $d_i \ge 2$ for all $i$, then $r<n$. If $r<n$, then there exists a (not necessarily) homogeneous
polynomial map $f:{\bf K}^r\rightarrow {\bf K}^r$ of degree
(at most) $d$ such that $f$ and $F$ are GZ-paired through some matrices $B\in M_{r,n}({\bf K})$ and
$C\in M_{n,r}({\bf K})$.
\end{cor}
\begin{proof}
The claim that $r < n$ under the given conditions follows by looking at the leading homogeneous part of
$\det JF \in {\bf K}$. So assume that the condition $r<n$ of the second claim is fulfilled.

By Lemma \ref{3.2}, $r = n-\dim_{\bf K} (\ker_{\bf K} JH)$.
Hence by the preceding theorem, there exists a (not necessarily) homogeneous polynomial map $f$ of degree
(at most) $d$ in dimension $r$, such that $f$ and $F$ are weakly GZ-paired through some
matrices $B\in M_{r,n}({\bf K})$ and $C\in M_{n,r}({\bf K})$. Since $\ker B\subseteq \ker A$ and
$\rank B=r=\rank A$, we have $\ker B=\ker A$. Thus $f$ and $F$ are GZ-paired.
\end{proof}

\begin{lem} \label{3.14}
Suppose that $f$ and $F=X+H$ are weakly GZ-paired through matrices $B\in M_{r,n}({\bf K})$ and
$C\in M_{n,r}({\bf K})$. Then we have the following.
\begin{enumerate}
\item[i)] If $y\in {\bf K}^n$ and $y_0 \in \ker_{\bf K} JH$, then $F(y+y_0)=F(y)+y_0$.
\item[ii)] If $y\in{\bf K}^n$, then $CBy-y\in \ker_{\bf K} B\subseteq \ker_{\bf K} JH$.
\end{enumerate}
\end{lem}
\begin{proof}
$ii)$ follows from $B(CBy-y)=BCBy-By=I_rBy-By=0$. To prove $i)$, notice first that $F(y+y_0)=F(y)+y_0$
is equivalent to $H(y+y_0)=H(y)$. Next we have
$$
H_i(y+y_0)-H_i(y)=\int_0^1 \Big(\frac{d}{dt} H_i(y+t y_0)\Big)dt
$$
and
$$
\frac{d}{dt} H_i(y+t y_0)= JH_i(y+t y_0)\cdot y_0=0
$$
because $y_0\in \ker_{\bf K} JH$. So $H_i(y+y_0)=H_i(y)$ for all $i$, which gives $i)$.
\end{proof}

\begin{cor} \label{our6.4.4}
Suppose that $f$ and $F=X+H$ are weakly GZ-paired through matrices $B\in M_{r,n}({\bf K})$ and
$C\in M_{n,r}({\bf K})$. Then
\begin{enumerate}
\item[i)] $BF(CBX) = BF(X)$, and
\item[ii)] $H(CBX) = H(X)$.
\end{enumerate}
\end{cor}

\begin{proof}
Using Lemma \ref{3.14}, we obtain that for all $y \in {\bf K}^n$,
$$
F(CBy) = F(y + (CBy - y)) = F(y) + (CBy - y)
$$
so $BF(CBX) = BF(X) + B(CBX - X) = BF(X)$ and $H(CBX) = F(CBX) - CBX =
F(X) + (CBX - X) - CBX = H(X)$.
\end{proof}

\begin{prop} \label{our6.4.8}
Suppose that $f$ and $F=X+H$ are (weakly) GZ-paired through matrices $B\in M_{r,n}({\bf K})$ and
$C\in M_{n,r}({\bf K})$. If either $f$ or $F$ is invertible, then they are both invertible. Furthermore,
$f^{-1}$ and $F^{-1}$ are (weakly) GZ-paired through $B$ and $C$ as well, and
$F^{-1}=X-H(Cf^{-1}(BX))$.
\end{prop}
\begin{proof}
We first show that $F$ is invertible if $f$ is invertible.
Hence suppose that $f$ is invertible. We show that $F^{*} := X-H(Cf^{-1}(BX))$ is the inverse of $F$.

By using $BC = I_r$, Corollary \ref{our6.4.4}, and $f = BF(Cx)$,
we can deduce that
\begin{align*}
F^{*}(F) &= F-H\big(Cf^{-1}(BF(X))\big) = F-H\big(Cf^{-1}(BF(CBX))\big) \\
&= F-H\big(C f^{-1}(f(BX))\big) = F-H(CBX) = F-H = X
\end{align*}
So $F^{*}$ is a left inverse of $F$. Furthermore, it follows from $BC = I_r$ that
\begin{align*}
B F^{*} (Cx) &= BC x - B H(C f^{-1} (BCx)) = x - B\big(F(C f^{-1}(x)) - C f^{-1}(x)\big) \\
&= x - BF(Cf^{-1}(x)) + BC f^{-1}(x) = x - f(f^{-1}(x)) + f^{-1}(x) = f^{-1}(x)
\end{align*}
Since additionally $B$ is a right factor of $J(-H(Cf^{-1}(BX))) = J(F^{*}-X)$, it follows
that $f^{-1}$ and $F^{*}$ are weakly GZ-paired through $B$ and $C$ as well.
So $F^{*}$ has a left inverse $F^{**}$, and $F^{**}(X) = F^{**}(F^{*}(F(X))) = F$
by associativity of composition, i.e.\@ $F^{-1}$ exists and is equal to $F^{*}$.

We next show that $f$ is invertible if $F$ is invertible.
Hence suppose that $F$ is invertible. We show that $f^{*} := B F^{-1} (Cx)$ is the inverse of $f$.

By using $f = BF(Cx)$, $i)$ of Corollary \ref{our6.4.4}, and $BC = I_r$, we can deduce that
$$
f(f^{*}(x)) = B F(CB F^{-1}(Cx)) = B F(F^{-1}(Cx)) = BCx = x
$$
So $f$ has a right inverse $f^{*}$. Since
\begin{equation} \label{fstareq}
\begin{split}
\ker B &= \ker (B|_{F^{-1}}) \subseteq \ker \big((JH)|_{F^{-1}}\big)
= \ker \big((JF)|_{F^{-1}} - I_n\big) \\
&= \ker \big(I_n - (JF)|^{-1}_{F^{-1}}\big) = \ker J(X-F^{-1})
\end{split}
\end{equation}
it follows that $f^{*}$ and $F^{-1}$ are weakly GZ-paired through $B$ and $C$ as well.
So $f^{*}$ has a right inverse $f^{**}$, and $f^{**}(x) = f(f^{*}(f^{**}(x))) = f(x)$
by associativity of composition, i.e.\@ $f^{-1}$ exists and is equal to $f^{*}$.

Finally, it remains to show that $f^{-1}$ and $F^{-1}$ are (weakly) GZ-paired through $B$ and $C$
as well. We have already deduced above that $f^{-1}$ and $F^{-1}$ are weakly GZ-paired
through $B$ and $C$. If $f$ and $F$ are GZ-paired through $B$ and $C$, then the inclusion in
\eqref{fstareq} is an equality. So $f^{-1}$ and $F^{-1}$ are GZ-paired as well in this case,
which completes the proof.
\end{proof}

\begin{rem}
Lemma \ref{3.14} is somewhat similar to \cite{MR1790619} (Lemma 6.4.4), which is valid with weakly
GZ-pairing as well. The purpose of Lemma \ref{3.14} and Corollary \ref{our6.4.4} is to replace
\cite{MR1790619} (Lemma 6.4.4) if not all components of $H$ are powers of a linear form.
We used this replacement in the proof of Proposition \ref{our6.4.8}, which in turn replaces the
cases where $F$ and $f$ are invertible, of \cite{MR1790619} (Proposition 6.4.7 iii)) and
\cite{MR1790619} (Proposition 6.4.8 iii)) respectively.
\end{rem}

\section{The JC in dimension \mathversion{bold}$r$}

The following theorem is a special case of (2) $\Rightarrow$ (3) of \cite{1203.6615} (Theorem 4.2), in which
$h$ may be any polynomial.

\begin{thm}
 Let $F=X+H$ such that $H_i \in {\bf K}[h]$ for all $i$, for a fixed linear form $h$.
If $\det JF=1$, then $F$ is invertible. More precisely, $F$ is linearly triangularizable.
\end{thm}
 \begin{proof}
Since $\ker_{\bf K} JH \supseteq \ker_{\bf K} Jh$ and $\ker Jh$ is generated by at least $n-1$ vectors over $K$,
we deduce that $\dim_{\bf K} (\ker_{\bf K} JH) \geq \dim_{\bf K} (\ker_{\bf K} Jh) =
\dim (\ker Jh)\ge n-1$.
So $1 \ge n - \dim_{\bf K} (\ker_{\bf K} JH)$, and it follows that with $F$, a polynomial mapping $f$ in
  dimension $1$ is weakly GZ-paired (see Theorem \ref{3.4}), say through the row matrix
$B$ and the column matrix $C$.

Since $C \notin \ker B$ and $\ker B$ is generated by $n-1$ vectors over $K$,
there exists an invertible matrix $T \in M_{n}({\bf K})$ such that the first column of $T$ is just $C$ and
all subsequent columns of $T$ are contained in $\ker B$. Since $\ker B \subseteq \ker JH$, only the first column
of $JH \cdot T$ might be nonzero. Using that
$$
J(H(TX)) = (JH)|_{TX} \cdot T = (JH \cdot T)|_{TX}
$$
we deduce that $J(T^{-1} F(TX)) = T^{-1} J(F(TX))$ is a lower triangular matrix with
ones on the diagonal, except maybe the leading diagonal entry. But that entry is also one
because $\det J(T^{-1} F(TX)) = (\det JF)|_{TX} = 1$. Hence $F$ is linearly triangularizable.
In particular, $F$ is a composition of $n-1$ elementary polynomial maps, so $F$ is invertible.
\end{proof}

\begin{rem}
Notice that in the above proof, the first row of $T^{-1}$ is just $B$.
Furthermore, $f$ and $F$ are GZ-paired, if and only if $\deg H \ge 1$.
This is because $\rank JH = 1 = \rank B$, if and only if $\deg H \ge 1$.
\end{rem}

\begin{prob}\label{3.7}
Let $F=X+H$ and $H_i=(A_iX)^{d_i}$ for $1\leq i\leq n$.
If $\det JF=1$ and $\rank (A)\leq 2$, then $F$ is invertible.
\end{prob}

\begin{thm} \label{3.8}
Problem \ref{3.7} is equivalent to the Jacobian Conjecture in dimension 2.
\end{thm}

\begin{proof}
We use the invertibility equivalence of GZ-pairing, which we proved in Proposition \ref{our6.4.8}.
From Corollary \ref{3.5}, we can deduce that with $F$, a polynomial map in dimension 2 is GZ-paired,
if $\rank (A)\leq 2$. So
if the Jacobian Conjecture is true in dimension 2, then Problem \ref{3.7} has an affirmative answer as well.

On the other hand, if $f:{\bf{K}}^2\rightarrow {\bf{K}}^2$ is a polynomial map, then
$f$ is GZ-paired with a non-homogeneous power-linear map $F_A=X+H$ such that $H_i = (A_iX)^{d_i}$
for each $i$, with
$\rank A\leq2$ (see Theorem \ref{3.3}). So if Problem \ref{3.7} has an affirmative answer,
then the Jacobian Conjecture is true in dimension 2.
\end{proof}

Next, we get a generalized statement of Problem \ref{3.7}.

\begin{prob} \label{3.9}
Let $F=X+H$ and $H_i=(a_{i1}x_1+a_{i2}x_2+\cdots+a_{in}x_n)^{d_i}$ for $1\leq i\leq n$.
If $\det JF=1$ and $\rank (A)\leq r$, then $F$ is invertible.
\end{prob}

\begin{rem}
We can assume $d_i\geq 2$ for all $1\leq i\leq n$ in Problems \ref{3.7} and \ref{3.9}. This is because we can
obtain $d_i\geq 2$ for all $1\leq i\leq n$ if we replace $F$ by $F(LX) - c$ for a suitable linear map $L$ and
a suitable $c \in {\bf K}^n$.
\end{rem}

\begin{thm}
Problem \ref{3.9} is equivalent to the Jacobian Conjecture in dimension $r$.
\end{thm}
\begin{proof}
Similar to the proof of Theorem \ref{3.8}
\end{proof}
Next, we give some partial results about Problem \ref{3.9}.

\begin{thm} \label{3.12}
Let $F=X+H$, where $H=(H_1,H_2,\dots,H_n)^t$ and $H_i$ is a homogeneous polynomial of
degree $d_i$ for $1\leq i\leq n$. If $\det (I+DJH|_{a})\ne 0$ for each
$\lambda \in \mathbf{K} \setminus \{1\}$ and every $a \in \mathbf{K}^n$, where
$$
D =\frac1{\lambda-1}
\diag \Big(\frac1{d_1}(\lambda^{d_1}-1),\frac1{d_2}(\lambda^{d_2}-1),\ldots,\frac1{d_n}(\lambda^{d_n}-1)\Big)
$$
then $F$ is injective on every line that passes through the origin.
More precisely, for each $\lambda \in \mathbf{K} \setminus \{1\}$ and every $a \in \mathbf{K}$,
we have $(I+DJH|_{a})a = 0$, if and only if $F(a) = F(\lambda a)$.

In particular, if $d_1=d_2=\cdots=d_n=d \ge 2$, then $\det JF=1$ is equivalent to
$\det (I+DJH)=1$, so homogeneous Keller maps are injective on lines through the origin.
\end{thm}

\begin{proof}
Since $H_i=d_i^{-1}\sum_{j=1}^{n}x_{j}H_{x_{j}}$ for $1\leq i\leq n$, we have $F=(I+D'JH)X$,
where $D'=\diag (d_1^{-1},d_2^{-1},\ldots,d_n^{-1})$.
Take $\lambda \in {\bf K} \setminus \{1\}$, $a\in {\bf K}^n$ and $b=\lambda a$.
Then $F(a) = F(b)$ is equivalent to
$$(I+D'JH_{a})a=(I+D'JH_{b})b,$$
which in turn is equivalent to
$$(I+D'JH_{a})a=(\lambda I+D'\Lambda JH_{a})a,$$
where $\Lambda=\diag (\lambda^{d_1},\lambda^{d_2},\cdots,\lambda^{d_n})$. That is,
$$[(\lambda-1)I+D'(\Lambda-I)JH_{a}]a=0.$$
Since $(\lambda-1)^{-1}D'(\Lambda-I)=D$, we see that $F(a) = F(b)$, if and only if
$$[I + DJH|_a]a = 0,$$ as desired.
\end{proof}

\begin{rem}
If $d_1=d_2=\cdots=d_n=d \ge 2$, then Theorem \ref{3.12} is similar to \cite{MR3133981}
(Proposition 1.1).
\end{rem}

\begin{lem} \label{3.16}
Suppose that $f$ and $F$ are weakly GZ-paired through matrices $B$ and $C$.
Let $a\in {\bf K}^n$ and $\lambda, \mu \in {\bf K}$. If $F(\lambda(a+b))\neq F(\mu(a+b))$
for all $b \in \ker_{\bf K} B$, then $f(\lambda Ba)\neq f(\mu Ba)$.
\end{lem}

\begin{proof}
Suppose that $f(\lambda Ba)= f(\mu Ba)$ and $\lambda \neq \mu$. By definition \ref{3.1},
we have $BF(\lambda CBa)=BF(\mu CBa)$. Hence $F(\lambda CBa)-F(\mu CBa) \in \ker_{\bf K} B$.
From $i)$ of Corollary \ref{our6.4.4}, it follows that $F(\lambda a) - F(\lambda CBa)
\in \ker_{\bf K} B$ and $F(\mu a) - F(\mu CBa) \in \ker_{\bf K} B$. Consequently,
$F(\lambda a) - F(\mu a) \in \ker_{\bf K} B$, say that
$$
F(\lambda a) - F(\mu a) = (\mu - \lambda) b
$$
where $b \in \ker_{\bf K} B$. By adding $F(\mu a) + \lambda b$ on both sides,
we get $F(\lambda a) + \lambda b = F(\mu a) + \mu b$. Hence we have
$F(\lambda(a+b)) = F(\mu(a+b))$ on account of $i)$ of Lemma \ref{3.14}.
This gives the desired result.
\end{proof}

\begin{thm} \label{3.17}
Let $F=X+H$, where $H=(H_1,H_2,\dots,H_n)^t$ and $H_i$ is a homogeneous polynomial of degree
$d_i$ for $1\leq i\leq n$. If $\ker_{\bf K} JH$ is a space of dimension $\geq n-2$,
and $\det (I+DJH) = 1$ for every $\lambda \in {\bf K} \setminus\{1\}$, where
$$
D =\frac1{\lambda-1}
\diag \Big(\frac1{d_1}(\lambda^{d_1}-1),\frac1{d_2}(\lambda^{d_2}-1),\ldots,\frac1{d_n}(\lambda^{d_n}-1)\Big)
$$
then $F$ is invertible.
\end{thm}

\begin{proof}
Assume first that ${\bf K}$ is algebraically closed.
Since $\ker_{\bf K} JH$ is a space of dimension $\geq n-2$, it follows that with
$F$, a polynomial map $f$ in dimension 2 is weakly GZ-paired, say through
matrices $B$ and $C$. Suppose that $f$ is not a Keller map.
Then we can find $a, b \in {\bf K}^2$, of which $b$ is nonzero, such that $b^t (Jf)|_{a} B$
is the zero matrix. By $BC = I_2$ and $i)$ of Corollary \ref{our6.4.4}
\begin{align*}
b^t \cdot (Jf)|_{a} \cdot B &= b^t \cdot (Jf)|_{BCa} \cdot B = b^t \cdot \big(J(f(BX))\big)\big|_{Ca} \\
&= b^t \cdot \big(J(BF(CBX))\big)\big|_{Ca} = b^t \cdot \big(J(BF)\big)\big|_{Ca}
= b^t B \cdot (JF)|_{Ca}
\end{align*}
so $b^t B (JF)|_{Ca}$ is the zero matrix as well. It follows that $F$ is not
a Keller map. But $\det JF=\det (I+IJH)= \lim_{\lambda\rightarrow 1} \det (I+DJH) = 1$.
Hence both $f$ and $F$ are Keller maps.

By Theorem \ref{3.12}, $F$ is injective on the lines that pass through the origin. Using
$\rank B = 2$,
Lemma \ref{3.16} subsequently gives that $f$ is injective on the lines that pass through the origin.
Thus $f$ is invertible, see \cite{MR1280980}. By the invertibility equivalence of GZ-pairing,
which we proved in Proposition \ref{our6.4.8}, $F$ is invertible.

Assume next that ${\bf K}$ is not algebraically closed, and let ${\bf \bar{K}}$
be the algebraic closure of ${\bf K}$. Since ${\bf K}$ is infinite, the identity
$\det (I+DJH) = 1$ for all $\lambda \ne 1$ remains
valid if we replace $\lambda$ by an indeterminate (so that $D$ becomes a rational function).
Hence $\det (I+DJH) = 1$ for all $\lambda \in {\bf \bar{K}}$.
Furthermore, $\ker_{\bf \bar{K}} JH$ is a space of
dimension $\geq n-2$, because it is generated by the same vectors over ${\bf K}$ as
$\ker_{\bf K} JH$. So $F$ is invertible over ${\bf \bar{K}}$. On account of
\cite{MR1790619} (Proposition 1.1.1), $F$ is invertible over ${\bf K}$ as well.
\end{proof}

\begin{cor}
Let $F=X+H$, where $H=(H_1,H_2,\dots,H_n)^t$ and $H_i=(a_{i1}x_1+a_{i2}x_2+\cdots+a_{in}x_n)^{d_i}$
for $1\leq i\leq n$. If $\det (I+DJH)=1$, where
$$
D =\frac1{\lambda-1}
\diag \Big(\frac1{d_1}(\lambda^{d_1}-1),\frac1{d_2}(\lambda^{d_2}-1),\ldots,\frac1{d_n}(\lambda^{d_n}-1)\Big)
$$
for any $\lambda\neq 1$, $\lambda \in {\bf K}$ and $\rank A\leq 2$, then $F$ is invertible.
\end{cor}

\begin{proof}
Since $\rank A\leq 2$, the space $\ker_K JH$ has dimension
$\geq n-2$. Hence the conclusion follows from Theorem \ref{3.17}.
\end{proof}

\section{A bound for the degree of the inverse of some special polynomial maps}

We start with a proposition that gives a connection between weak GZ-pairing and
the degree of the inverse.

\begin{prop}
Suppose that $f$ and $F$ are weakly GZ-paired and have degree at most $d$.
If either $f$ or $F$ is invertible, then they are both invertible and
$\deg f^{-1} \le \deg (F^{-1})\leq d\cdot\deg f^{-1}$.
\end{prop}
\begin{proof}
Say that $f$ and $F$ are weakly GZ-paired through matrices $B\in M_{r,n}({\bf K})$ and
$C\in M_{n,r}({\bf K})$. Suppose that either $f$ or $F$ is invertible. Then Proposition \ref{our6.4.8},
tells us that they are both invertible and that $F^{-1}=X-H(Cf^{-1}(BX))$, where $H(X)=F(X)-X$.
Since $\deg H(X)\leq d$, we have
$\deg (F^{-1})\leq d\cdot\deg f^{-1}$ indeed.
\end{proof}

The following theorem is a generalization of \cite{MR2419134} (Theorem 3(1)) to the case where $H_i$
is not a power of a linear form for some $i$.

\begin{thm} \label{idxth}
Assume that $f = x + h$ and $F = X + H$ are polynomial maps in dimensions $r < n$ respectively.
Then we have the following.
\begin{enumerate}
\item[i)] If $f$ and $F$ are weakly GZ-paired, then $Jh^k = 0$ implies $JH^{k+1} = 0$.
\item[ii)] If $f$ and $F$ are GZ-paired, then $Jh^k = 0$, if and only if $JH^{k+1} = 0$.
\end{enumerate}
\end{thm}

\begin{proof}
Assume that $f$ and $F$ are weakly GZ-paired through matrices $B$ and $C$.
Since $x + h = B\big(Cx+H(Cx)\big)$ and $BC=I_r$, we have $h = BH(Cx)$. Hence $Jh = B JH|_{Cx} C$
and $JH|_{Cx} \cdot C Jh^k B = (JH|_{Cx} \cdot CB)^{k+1}$. If we substitute $x = BX$
and use $iii)$ of Lemma \ref{3.14}, then we obtain
\begin{equation} \label{idxeq}
\begin{split}
(JH)|_{CBX} \cdot C (Jh)|_{BX}^k B &= \big((JH)|_{CBX} \cdot CB\big)^{k+1}\\
&= J(H(CBX))^{k+1} = JH^{k+1}
\end{split}
\end{equation}
So if $Jh^k = 0$, then also $(Jh)|_{BX}^k = 0$ and $JH^{k+1} = 0$, which gives i).

Suppose next that $f$ and $F$ are GZ-paired and $JH^{k+1} = 0$.
From \eqref{idxeq}, we obtain that $(JH)|_{CBX} \cdot C (Jh)|_{BX}^k B = 0$, i.e.\@
\begin{equation} \label{kereq}
(Jh)|_{BX}^k B \in \ker \big((JH)|_{CBX} \cdot C\big)
\end{equation}
Using $BC = I_r$, $iii)$ of Lemma \ref{3.14} and $\ker JH = \ker B$, we see that
\begin{align*}
\ker \big((JH)|_{CBX} \cdot C\big) &= \ker \big((JH)|_{CBX} \cdot CBC\big)
= \ker \big(J(H(CBX)) \cdot C\big) \\
&= \ker (JH \cdot C) = \ker (B \cdot C) = \ker I_r
\end{align*}
Since $\ker I_r$ is trivial, we deduce from \eqref{kereq} that $Jh|_{BX}^k B = 0B$.
Hence by $B \cdot C = I_r$, $Jh^k = (Jh)|_{B \cdot Cx} B \cdot C = (0B)|_{Cx} \cdot C = 0$,
as desired.
\end{proof}

\begin{thm} \label{4.1}
Let $F=X+H$ be a polynomial map, such that $H$ is homogeneous.
\begin{enumerate}

\item[i)] If $JH^2 = 0$, then $F$ is invertible and $F^{-1}=X-H$.

\item[ii)] If $\dim (\ker JH) = \dim_{\bf K} (\ker_{\bf K} JH)$ and $JH^3=0$, then $F$ is invertible
and $F^{-1}=X-H(X-H)$.
\end{enumerate}
\end{thm}

\begin{proof}
Write $d := \deg H$.
\begin{enumerate}

\item[i)] Suppose that $JH^2 = 0$. We shall show that $JH \cdot H = 0$.
If $d = 0$ then $JH = 0$. If $d \ge 1$, then by Euler's homogeneous function theorem,
$JH \cdot H = d^{-1} JH^2 X = 0$. So $JH \cdot H = 0$ and by \cite{homokema}
(Proposition 3.1.2), $2X - F = X - H$ is the inverse polynomial map of $F$.

\item[ii)] Suppose that $JH^3=0$. Then we have $r := n - \dim_{\bf K} (\ker_{\bf K} JH)
= n - \dim (\ker JH) = \rank JH < n$. Hence by Theorem \ref{3.4},
there exists a polynomial map $f$ such that $f$ and $F$ are weakly
GZ-paired through some matrices $B\in M_{r,n}({\bf K})$ and $C\in M_{n,r}({\bf K})$.
Since $\rank B \le r \le \rank JH$ and $\ker B \subseteq \ker JH$,
we have $\ker B = \ker JH$. So $f$ and $F$ are GZ-paired through $B$ and $C$.

Write $h = BH(CX)$. Then $h$ is homogeneous of the same degree $d$ as $H$ is.
From Theorem \ref{idxth}, $Jh^2 = 0$ follows. Hence by i), $2x - f = x - h$
is the inverse polynomial map of $f$, where $x = (x_1, x_2, \ldots, x_r)^t$.
Using Proposition \ref{our6.4.8}, $f^{-1} = 2x - f$, $i)$ of definition \ref{3.1}
and $F = X + H$, in that order, we obtain
\begin{align*}
F^{-1} &= X - H\big(C f^{-1}(BX)\big) = X - H\big(C(2BX - f(BX))\big) \\
       &= X - H\big(C (2BX - BF(CBX))\big) = X - H\big(CBX - CBH(CBX)\big)
\end{align*}
Now $F^{-1}=X-H(X-H)$ follows by applying $ii)$ of Corollary \ref{our6.4.4} twice on the right hand side.
\qedhere

\end{enumerate}
\end{proof}

\begin{cor} \label{4.2}
Let $F=X+H$ be a homogeneous power linear map. If $JH^k=0$ for some $k \le 3$,
then $F$ is invertible and $F^{-1}=X-H(X-H)$. Furthermore, $\deg (F^{-1})\leq (\deg F)^k$.
\end{cor}

\begin{proof} By power linearity, $\ker JH = \ker A$ for some $A \in M_n({\bf K})$.
Hence the previous theorem gives the desired result.
\end{proof}

\begin{ques} \label{4.3}
Let $F=X+H$ be a homogeneous power linear map. If $F$ is invertible and $JH^{k+1}=0$, then
$\deg (F^{-1})\leq (\deg F)^k$.
\end{ques}

We see from Corollary \ref{4.2} that Question \ref{4.3} has an affirmative answer for $k\leq 2$.
However, it is not true for $k\geq 3$. We will give a counterexample below.

\begin{thm} \label{4.4}
Suppose that $f$ and $F$ are (weakly) GZ-paired through matrices $B$ and $C$.
Then $f$ and $\tilde{F}:= (F, x_{n+1} + (B_i X)^d)$ are (weakly) GZ-paired as well. Furthermore,
if $f$ is invertible, then $\tilde{F}$ is also invertible, and the degree of the last component of
$\tilde{F}^{-1}$ is $d$ times the degree of the $i$-th component of $f^{-1}$.
\end{thm}
\begin{proof}
 One can easily see that $f$ and $\tilde{F}$ are (weakly) GZ-paired through matrices $\tilde{B}$ and
$\tilde{C}$, where $\tilde{B}$ is obtained from $B$ by adding a zero column on the right, and $\tilde{C}$ is
obtained from $C$ by adding an arbitrary row on the bottom. By Proposition \ref{our6.4.8} and by definition of
$\tilde{B}$ and $\tilde{C}$, the last component of $\tilde{F}^{-1}(X)$ equals
$x_{n+1} - \big(B_i C f^{-1}(BX)\big)^d$, which by $BC=I_r$ simplifies to $x_{n+1} - (f^{-1}(BX))_i^d$.
The degree of $(f^{-1}(BX))_i$ is equal to that $(f^{-1}(x))_i$, because
$$
  \deg \big(f^{-1}(x)\big)_i = \deg \big(f^{-1}(BCx)\big)_i
    \leq \deg \big(f^{-1}(BX)\big)_i
    \leq \deg \big(f^{-1}(x)\big)_i
$$
 which completes the proof.
\end{proof}

Next we give a counterexample of Question \ref{4.3}.
\begin{exa}[Furter] \label{4.5}
Let $x$ and $h$ be given by $x = (x_1,x_2, \ldots, x_6)^t$ and
$h=(2x_2x_6-2x_3^2-x_4x_5, 2x_3x_5-x_4x_6, x_5x_6, x_5^2, x_6^2, 0)^t$, and $f=x+h$.
\end{exa}

Following Theorem \ref{3.3}, we get a homogeneous power linear $F$, say in dimension $n$,
 with which $f$ is GZ-paired. Suppose $f$ and $F$ are weakly
GZ-paired through matrices $B$ and $C$. By Theorem \ref{4.4}, we get
$\tilde{F}=(F,x_{n+1}+(B_1X)^2)$ and $\deg(\tilde{F})^{-1}\geq 2\deg (f^{-1})_1$.
Write $\tilde{F}=\tilde{X}+\tilde{H}$. Since $Jh^3=0$, we have $J\tilde{H}^4=0$ by Theorem \ref{idxth}.
It is easy to compute that $\deg(f^{-1})_1=6$. Thus $\deg(\tilde{F}^{-1})\geq 12>8=2^{4-1}$.

\begin{rem}
Example \ref{4.5} also shows that the assumption that $\dim (\ker JH) = \dim (\ker_{\bf K} JH)$ in
Theorem \ref{4.1} is necessary.
\end{rem}

Notice that $\rank JH\leq k$ implies $JH^{k+1}=0$ when $JH$ is nilpotent.
We get the following question if we replace $JH^{k+1}=0$ by $\rank JH\leq k$ in Question
\ref{4.3}.

\begin{ques} \label{4.7}
Let $F=X+H$ be a polynomial map over ${\bf K}[x_1,x_2,\ldots,x_n]$. If $F$ is invertible and
$\rank JH\leq k$, then $\deg (F^{-1})\leq (\deg F)^k$.
\end{ques}

In \cite{MR3133981} (Theorem 3.4), we showed that Question \ref{4.7} has an affirmative answer if
$\ker JH = \ker_{\bf K} JH$, because $\dim (\ker JH) = n-\rank JH \geq n-k$.
This is in particular the case when $F$ is (non-homogeneous) power linear, see \cite{MR3133981} (Theorem 3.5).

\begin{thm}
If $\rank JH \leq 1$ or $\rank JH \geq n-1$,
then Question \ref{4.7} has an affirmative answer.
\end{thm}

\begin{proof}
In \cite{MR0663785}, it has been proved that the degree of the inverse of any invertible polynomial map $F$
is at most $(\deg F)^{n-1}$. This gives the case $\rank JH \geq n-1$.

So assume that $\rank JH \leq 1$. Reading the proof of \cite{1203.6615} (Theorem 4.2),
we see that there exists a $T \in \operatorname{GL}_n(K)$ such that $T^{-1}H(TX)$ is of the form
$$
\big[c_1,c_2,\ldots,c_s,\lambda x_{s+1} + g, h_{s+2}(\lambda x_{s+1} + g), h_{s+3}(\lambda x_{s+1} + g),
\ldots, h_n(\lambda x_{s+1} + g)\big]^t,
$$
where $0 \le s \le n-1$, $c_i \in \mathbf{K}$ for all $i$, $\lambda \in \mathbf{K} \setminus \{-1\}$
and $g \in K[x_1, x_2, \ldots,x_s]$. One can verify that
$$
\Big(\frac1{\lambda+1}(\lambda x_{s+1} + \tilde{g})\Big)\Big|_{F} = \lambda x_{s+1} + g
$$
where $\tilde{g} = g(x_1-c_1, x_2-c_2, \ldots, x_s - c_s)$, and that the inverse of $T^{-1} F(TX)$ is
\begin{align*}
&\Big[x_1-c_1,x_2-c_2,\ldots,x_s-c_s,x_{s+1}-\Big(\frac1{\lambda+1}(\lambda x_{s+1} + \tilde{g})\Big),\\
&~x_{s+2}-h_{s+2}\Big(\frac1{\lambda+1}(\lambda x_{s+1} + \tilde{g})\Big),\ldots,
x_n-h_n\Big(\frac1{\lambda+1}(\lambda x_{s+1} + \tilde{g})\Big)\Big]^t.
\end{align*}
Hence $F^{-1}$ has the same degree as $F$ itself.
\end{proof}

\def\polhk#1{\setbox0=\hbox{#1}{\ooalign{\hidewidth
  \lower1.5ex\hbox{`}\hidewidth\crcr\unhbox0}}}

\end{document}